\documentclass[12pt,leqno]{article}
\usepackage{amssymb}

\newcommand{\Abar}{{\backslash\kern-8pt A}}

\title{\sf\LARGE BRINGING ERRORS INTO FOCUS}
\author{\sc Nicolas Bouleau\footnote{Ecole des Ponts, ENPC, ParisTech,  28 rue des Saints P\`eres, 75007 Paris; 
e-mail : {\tt bouleau@enpc.fr}}}
\date{\it Lecture at the colloquium of mathematics of the University of Strasbourg, France\\
 7th April 2006}

\begin{document}
\maketitle
{\Large

\noindent {\it Summary}\\

{\sf I. The dichotomy of small errors, {\large p 2},\\

$\!\!$II. Intrinsic error calculi, {\large p 4},\\

$\!\!\!\!$III. Complete and symmetric error calculus, {\large p 13}, \\

$\!\!\!\!$IV. The four bias operators, {\large p 21}, \\

$\!\!\!$V. Statistics and errors, {\large p 26}.}}\\

\newpage

{\bf\noindent\textsf{\LARGE I. The dichotomy of small errors.}}\\

{\large There are two kinds of small errors which do not follow the same differential calculus. \\

In order to understand this phenomenon, let us consider two researchers, two applied mathematicians say, who  endeavor to perform random simulation rigorously. Both are expert in using inversion or rejection methods so that they are able to simulate any probability law as soon as they can draw a real number randomly in the interval  $[0,1]$.

For this, congruence methods are of course available, statistically excellent, but instead they want to be able to assess the committed error :

- the first one draws the binary digits by heads or tails. 

- the second one uses  Polya's urn.\\

Let us look at the biases and the variances of the errors after  $n$ drawings.\\

\underline{In the first case} 
$\;x=0,a_1a_2a_3\cdots $ and we approach
$$x=\sum_{k=1}^\infty\frac{a_k}{2^k}\qquad{\mbox{by}}\qquad x_n=\sum_{k=1}^n\frac{a_k}{2^k}\qquad\qquad a_k\in\{0,1\}$$
The $a_i$ are independent random variables  with value 1 or 0 with probability $\frac{1}{2}$, 

\noindent$\mathcal{F}_n=\sigma(a_1,\ldots,a_n)$. 
Let us put 
$$\begin{array}{rl}
b_n&=\mathbb{E}[(x-x_n)|\mathcal{F}_n]\\
d_n&=\mathbb{E}[(x-x_n)^2|\mathcal{F}_n]\\
v_n&=d_n-b_n^2.
\end{array}
$$
We have
$$b_n=\sum_{n+1}^\infty\frac{1/2}{2^k}=\frac{1}{2^{n+1}},\qquad d_n=\frac{1}{3}\frac{1}{4^n},\qquad v_n=\frac{1}{3}\frac{1}{4^{n+1}}.$$

\indent\underline{In the second case},  let us recall the rule of Polya's urn.  At the beginning there are one white and one black ball in the urn and each time a ball is drawn out, it is put back into the urn together with another one of the same color. 

After $n$ drawings, the ratio of white balls in the urn may be written 

$$X_{n+1}(n+3)=X_n(n+2) +1_{\{U_{n+1}\leq X_n\}}$$
where $U_{n+1}$ is uniformly distributed on $[0,1]$ independent of $\mathcal{F}_n=\sigma(X_0,\ldots,X_n)$.

In other words 
$$X_{n+1}=X_n+\frac{1}{n+3}(1_{\{U_{n+1}\leq X_n\}}-X_n).$$
We see that $X_n$ is a bounded martingale which converges a.s. and in $L^p$, $p\in[1,\infty[$, to $X_\infty$ and it is not difficult to prove that when the initial composition of the urn is one white and one black ball,   $X_\infty$ is uniformly distributed on $[0,1]$.

We have for the bias and the variance

$$b_n=\mathbb{E}[X_\infty-X_n|\mathcal{F}_n]=0$$
$$v_n=d_n=\mathbb{E}[(X_\infty-X_n)^2|\mathcal{F}_n]\quad\mbox{ with }\quad\mathbb{E}[v_n]=\frac{1}{6n}+o(1/n).$$

We see that in the first case the variances are smaller than the biases, and in the second case the biases are smaller than the variances. \\

How  will this propagate trough the computations of our two modelers ?
\\

Let us write the Taylor expansion of a 
$\mathcal{C}^3$-function with bounded derivatives.
$$f(X)-f(X_n)=(X-X_n)f^\prime(X_n)+\frac{1}{2}(X-X_n)^2f^{\prime\prime}(X_n)+\frac{1}{6}(X-X_n)^3f^{\prime\prime\prime}(X_n+\theta(X-X_n))$$
$$\begin{array}{rcl}
B_n=&\mathbb{E}[f(X_\infty)-f(X_n)|\mathcal{F}_n]&=b_nf^\prime(X_n)+\frac{1}{2}d_nf^{\prime\prime}(X_n)+o(d_n)\\
D_n=&\mathbb{E}[(f(X_\infty)-f(X_n))^2|\mathcal{F}_n]&=d_n f^{\prime 2}(X_n)+o(d_n)
\end{array}$$
Three cases appear
 (supposing that the derivatives $f^\prime$ and $f^{\prime\prime}$ do not vanish for simplicity, or that the sets of their zeros are negligible for the law of $X_n$)~:

$1^o\!/$ When the variance is negligible with respect to the bias, $v_n\ll b_n$, (first researcher), we have also $d_n\ll b_n$, asymptotically the dominating term for the bias $B_n$ is the first term  $b_nf^\prime(X_n)$. We see that  $D_n$ is negligible before $B_n$ so that $V_n\ll B_n$,  the situation will endure. The only useful formula is 
$$\mathbb{E}[f(X_\infty)-f(X_n)|\mathcal{F}_n]=b_nf^\prime(X_n)+o(b_n).$$

$2^o\!/$ When the variance is of the same order of magnitude as the bias, this situation will persist ($D_n$ and $V_n$ are of the same order of magnitude as $B_n$).

$3^o\!/$ When the bias is negligible before the variance  $b_n\ll v_n$ (second researcher), the dominating term of the bias becomes $\frac{1}{2}d_nf^{\prime\prime}(X_n)$ which is equivalent to $\frac{1}{2}v_nf^{\prime\prime}(X_n)$ and we fall back into the second case ($2^o\!/$) where the bias  $B_n$ and the   variance $V_n$ are of the same order of magnitude.

In the cases $2^o\!/$ and $3^o\!/$,  $b_n^2$ may be neglected before $v_n$ and the useful formulae are
$$\begin{array}{rl}
\mathbb{E}[f(X_\infty)-f(X_n)|\mathcal{F}_n]&=b_nf^\prime(X_n)+\frac{1}{2}v_nf^{\prime\prime}(X_n)+o(v_n)\\
\mathbb{E}[(f(X_\infty)-f(X_n))^2|\mathcal{F}_n]&=v_n f^{\prime 2}(X_n)+o(v_n).
\end{array}$$

We see that our first modeler may content himself with an error calculus based on the first derivatives.

Instead, the second modeler, who uses Polya's urn, needs an error calculus for the bias and the variances. 

- the calculus for the variances is a first order calculus, 

- the one for the biases is a second order calculus which uses the calculus on variances. \\

In practice, generally, we do not control the nature of the errors. In a modelization,  \underline{errors on data are exogenous}, we do not know exactly from where they come. Therefore it is wise to do as if we were in the second case,  especially to take in account the stochastic nature of the errors and the non-linearity of the model. 

This may be said otherwise : to consider that a sensitivity analysis on the data or on the parameters may be done with the first derivatives, is to suppose that the errors on these quantities are representable by their mean. If, for more and more accurate measurements, the errors are thought as a stochastic process, this amount to consider it satisfies the ordinary differential calculus instead of the Ito calculus.  We shall come back to this analogy later on.  \\

{\bf\noindent\textsf{\LARGE II. Intrinsic error calculi.}}\\

In the error calculus on biases and variances, the calculus on variances do not involve the biases, it is natural to begin with it, and this takes us back to Gauss at the beginning of the   XIX-th century.\\

\noindent\underline{Gauss' error calculus on variances and covariances}.\\

Twelve years after his demonstration that, under some hypotheses, error follow the normal law, Gauss is interested in the propagation of errors  (\textsl{Theoria combinationis} 1821). Given a quantity  $U=F(X_1,X_2,\ldots)$ function of other quantities  $X_1,X_2,\ldots$, he consider the problem of computing the quadratic error on  $U$ knowing the quadratic errors $\sigma_1^2,\sigma_2^2,\ldots$ on $X_1,X_2,\ldots$, assuming these errors are small and independent. 

His answer is the following : 

\begin{equation}
\label{(1)}
\sigma_U^2=(\frac{\partial F}{\partial X_1})^2\sigma_1^2+(\frac{\partial F}{\partial X_2})^2\sigma_2^2+\cdots
\end{equation}
and for another function  $V=G(X_1,X_2,\ldots)$ he gives the covariance of the errors on $U$ and $V$:
\begin{equation}
\label{(2)}
cov_{UV}=\frac{\partial F}{\partial X_1}\frac{\partial G}{\partial X_1}\sigma_1^2+
\frac{\partial F}{\partial X_2}\frac{\partial G}{\partial X_2}\sigma_2^2+\cdots
\end{equation} Gauss didn't study the propagation of biases, which, as we shall see, is more delicate. 
\vspace{.5cm}

\noindent\underline{The erring ways of textbooks :  {\tt "}ugly{\tt "} formulae}.\\

In spite of the works of Gauss, during the whole XIX-th century and also in the XX-th century, have been taught in mathematical or physical textbooks very ambiguous formulae like 
\begin{equation}
\label{(3)}
\Delta U=|\frac{\partial F}{\partial X_1}|\Delta X_1+|\frac{\partial F}{\partial X_2}|\Delta X_2+\cdots
\end{equation}
as, for example, in the course \textsl{Math\'ematiques g\'en\'erales} (1947) of Vessiot and Montel, where the  $\Delta X_1, \Delta X_2, \Delta U$, are the absolute value of the estimated error. Sometimes (3) is justified  (J. Taylor \textsl{An Introduction to Error Analysis}, University Science Books 1992) by the inequality \begin{equation}
\label{(4)}
\sqrt{(\frac{\partial F}{\partial X_1})^2\sigma_1^2+\cdots+(\frac{\partial F}{\partial X_k})^2\sigma_k^2}\leq 
|\frac{\partial F}{\partial X_1}|\sigma_1+\cdots
+|\frac{\partial F}{\partial X_k}|\sigma_k.
\end{equation} which would make Gauss' formulae (1) and (2) useless. 
In addition to the fact that in (3) it is not at all clear what are  $\Delta X_1, \Delta X_2, \Delta U$,  these formulae are \textit{clumsy}.

It is not a question of taste or aesthetics, it is an affair of symbolism and concepts : with (3),  $|\Delta U|$ generally depends on the way the function $F$ is written. By composing two linear maps with values in  $\mathbb{R}^2$ we see that {\it the identity map increases the errors} ... impossible to work  properly in such conditions. \\

This does not happen in the calculus of Gauss. To see this, we may introduce the differential operator 
$$L=\frac{1}{2}\sigma_1^2\frac{\partial^2}{\partial X_1^2}+\frac{1}{2}\sigma_2^2\frac{\partial^2}{\partial X_2^2}+\cdots
$$ and remark that (1) writes
$$\sigma_U^2=LF^2-2FLF.$$
Then the coherence is related to the transport of a differential operator by a function.  If $u$ and $v$ are regular injective mappings and if we denote  $\theta_u L$ the operator $\varphi\mapsto L(\varphi\circ u)\circ u^{-1}$, we have $\theta_{v\circ u}L=\theta_v(\theta_uL)$. We shall make this more precise later on.\\

\noindent\underline{Geometrization}.\\

Errors on $X_1,X_2,\ldots$ may be supposed non-independent and may depends on the values of  $X_1,X_2,\ldots$~: we consider a field of positive symmetric matrices  $\sigma_{ij}(x_1,x_2,\ldots)$ on $\mathbb{R}^d$ representing the conditional variances and covariances of the errors given the values  $x_1,x_2,\ldots$ of $X_1,X_2,\ldots$, and the error calculus becomes 
\begin{equation}
\label{(5)}
\sigma_F^2=\sum_{ij}\frac{\partial F}{\partial X_i}(x_1,x_2, \ldots)\frac{\partial F}{\partial X_j}(x_1,x_2, \ldots)\sigma_{ij}(x_1,x_2, \ldots)
\end{equation}

In order to geometrize the error calculus, i.e. to find a language for the biases and the variances which depends on the only mathematical objects and not on their written form, we proceed in two steps. 

First we argue on  $\mathbb{R}^d$, setting about it in such a way that if  $$F=f\circ g=k\circ h$$ the error depends only on  $F$. Then we shall extend this when  $F$ takes its values in a manifold.\\

\noindent{\bf a)} Let $F$ be a regular  function,  $\mathcal{C}^\infty$ say, of two variables  $x$ and $y$ say. Let us consider the increment of  $F$ between $(x,y)$ and $(x+dx, y+dy)$. Here  $dx$ and $dy$ are arbitrary increments. 

Taylor formula gives
$$
\Delta F=F(x+dx,y+dy)-F(x,y)=P_1(dx,dy)+\cdots+\frac{1}{n!}P_n(dx,dy)+(|dx|+|dy|)^n o(|dx|+|dy|)$$
where
$$P_k(dx,dy)=\sum_{p+q=k}\frac{k!}{p!q!}\frac{\partial^j F}{\partial^p x\partial^q y}dx^pdy^q$$
is a homogeneous polynomial of degree  $k$ in $dx,dy$ with coefficients depending on  $x,y$ which we shall call \underline{the differential of order $k$ of $F$} denoted $d^k F$.

From the formula
$$\Delta F=dF+\frac{1}{2}d^2F+\cdots+\frac{1}{n!}d^nF+(|dx|+|dy|)^no(|dx|+|dy|)$$
we deduce the differential calculus of interest for us~: 
if $x=f(u,v,w)$ and $y=g(u,v,w)$ and if
$Z(u,v,w)=F(f(u,v,w),g(u,v,w))$ we have
\begin{equation}
\label{(6)}
\left\{
\begin{array}{rl}
dZ&=P_1(du,dv,dw)=\frac{\partial F}{\partial x}dx+\frac{\partial F}{\partial y}dy\\
d^2Z&=P_2(du,dv,dw)=d^2F\quad{(\mbox{\tt\small quadratic terms in }}dx,dy) \\
& \hspace{4cm}+\frac{\partial F}{\partial x}d^2x+\frac{\partial F}{\partial y}d^2y\\
& \hspace{4cm}{(\mbox{\tt\small linear terms in }}d^2x,d^2y,\;{\mbox{\tt\small  quadratic in }}du,dv,dw)\\
&\\
&=\frac{\partial^2 F}{\partial x^2}(dx)^2+2\frac{\partial^2 F}{\partial y\partial y}dxdy
+\frac{\partial^2 F}{\partial y^2}(dy)^2+
\frac{\partial F}{\partial x}d^2x+\frac{\partial F}{\partial y}d^2y
\end{array}
\right.
\end{equation}
what may be remembered by writing 
$$
\begin{array}{rl}
d^2Z=d(dF)&=d(\frac{\partial F}{\partial x}dx+\frac{\partial F}{\partial y}dy)\\
&\\
&=d(\frac{\partial F}{\partial x})dx+\frac{\partial F}{\partial x}d^2x
+d(\frac{\partial F}{\partial y})dy+\frac{\partial F}{\partial y}d^2y.
\end{array}
$$

\noindent{\bf b)} \underline{Randomization}. Formulae (6) are  identically  valid when the increments  $du,dv,dw$ are supposed to be random. 

Simply the Landau symbols $o(.)$ are now random and  we must be carefull to the integrability when expectations are taken. 

Thus, we may write with capitals as is usual in probability theory  : if $Z=F(X,Y)$
\begin{equation}
\label{(7)}
\left\{
\begin{array}{rl}
dZ=&\frac{\partial F}{\partial x}dX+\frac{\partial F}{\partial y}dY\\
d^2Z=& \frac{\partial^2 F}{\partial x^2}(dX)^2+2\frac{\partial^2 F}{\partial y\partial y}dXdY
+\frac{\partial^2 F}{\partial y^2}(dY)^2+
\frac{\partial F}{\partial x}d^2X+\frac{\partial F}{\partial y}d^2Y.
\end{array}
\right.
\end{equation}
If the functions $f,g$ and $F$ have bounded derivatives and if the random variables $dU,dV,dW$ are square integrable, we see that $dX,dY$ and $dZ$ are square integrable and $d^2X,d^2Y$ and $d^2Z$ are in  $L^1$.\\

We introduce now the bias and variance operators   $A$ and $\Gamma$ associated to the erroneous random variable  $(X,Y)$ ~:\begin{equation}
\label{(8)}
\left\{
\begin{array}{rl}
A[F](x,y)=&\mathbb{E}[dZ+\!\frac{1}{2}d^2Z\;|\;X\!=\!x,Y\!=\!y]\\
&\\
\Gamma[F](x,y)=&\mathbb{E}[(dZ)^2\;|\;X\!=\!x,Y\!=\!y]
\end{array}
\right.
\end{equation}
Let us point out at once that in the formula giving  $A[F]$ we may or not put  $dZ$, we shall come back on this, we let it for the moment. \\

- Let us remark that if we apply $A$ to the function $F^2$, we obtain $d^2(F^2)=d(2FdF)=2(dF)^2+2Fd^2F$, so  that
\begin{equation}
\label{(9)}
\begin{array}{c}
A[F^2]=\mathbb{E}^{x,y}[2FdF+(dF)^2+Fd^2F]=2FA[F]+\Gamma[F]\\
\Gamma[F]=A[F^2]-2FA[F].
\end{array}
\end{equation}

- If we consider
$$H(X,Y)=\Phi(F_1(X,Y),F_2(X,Y),\ldots,F_d(X,Y))$$
we see by the change of variable formulae  (7) that
\begin{equation}
\label{(10)}
\left\{
\begin{array}{rl}
A[H]=&{\displaystyle\sum_{i=1}^d\frac{\partial \Phi}{\partial F_i}A[F_i]+\frac{1}{2}\sum_{i,j=1}^d\frac{\partial^2\Phi}{\partial F_i\partial F_j}\Gamma[F_i,F_j]}\\
&\\
\Gamma[H]=&{\displaystyle\sum_{i,j=1}^d\frac{\partial \Phi}{\partial F_j}\frac{\partial \Phi}{\partial F_j}\Gamma[F_i,F_j]}
\end{array}
\right.
\end{equation}
the second relation generalizes that of Gauss to the case of possibly correlated errors.\\

\noindent{\bf c)} \underline{The analogy of formulae  (10) with Ito's formula} for semi-martingales and their brackets is striking. It may be made more precise as follows~: 

Let us consider the vector field 
$$
\begin{array}{rcl}
A[X](x,y)& = & \mathbb{E}[dX+\frac{1}{2}d^2X|X=x,Y=y] \\
A[Y](x,y)& = & \mathbb{E}[dY+\frac{1}{2}d^2Y|X=x,Y=y]
\end{array}
$$
and the field of matrices
$$
\underline{\underline{\Gamma}}(x,y)=\left(\begin{array}{cc} \Gamma[X](x,y) & \Gamma[X,Y](x,y) \\\Gamma[X,Y](x,y) & \Gamma[Y](x,y)   \end{array}\right)
$$
and let $\Sigma(x,y)$ be a regular square root of $\underline{\underline{\Gamma}}$, let us consider the diffusion process $S$ solution to the stochastic differential equation
$$
\left\{
\begin{array}{rl}
dS_t^1=&\Sigma_{11}(S_t)dB^1_t+\Sigma_{12}(S_t)dB_t^2+A[X](S_t)dt\\
dS_t^2=&\Sigma_{21}(S_t)dB^1_t+\Sigma_{22}(S_t)dB_t^2+A[Y](S_t)dt
\end{array}\right.
$$
where $B=(B^1,B^2)$ is a standard Brownian motion. Ito's calculus gives 
$$
\begin{array}{rl}
dF_i(S)_t=&\frac{\partial F_i}{\partial X}(S_t)dS_t^1+\frac{\partial F_i}{\partial Y}(S_t)dS_t^2\\
&+\frac{1}{2}[\frac{\partial^2 F_i}{\partial X^2}(S_t)\Gamma[X](S_t)
+2\frac{\partial^2 F_i}{\partial X\partial Y}(S_t)\Gamma[X,Y](S_t)
+\frac{\partial^2 F_i}{\partial Y^2}(S_t)\Gamma[Y](S_t)]dt
\end{array}
$$
and also, denoting $(.)^\ast$ the continuous finite variation part \begin{equation}
\label{(11)}
\left\{
\begin{array}{l}
(dF_i(S)_t)^\ast=\frac{\partial F_i}{\partial X}(S_t)A[X](S_t)dt+\frac{\partial F_i}{\partial Y}(S_t)A[Y](S_t)dt\\
\qquad\qquad+\frac{1}{2}[\frac{\partial^2 F_i}{\partial X^2}(S_t)\Gamma[X](S_t)
+2\frac{\partial^2 F_i}{\partial X\partial Y}(S_t)\Gamma[X,Y](S_t)
+\frac{\partial^2 F_i}{\partial Y^2}(S_t)\Gamma[Y](S_t)]dt\\
{\mbox{and}}\\
d<F_i(S),F_j(S)>_t=\frac{\partial F_i}{\partial X}(S_t)\frac{\partial F_j}{\partial X}(S_t)\Gamma[X](S_t)dt + 
(\frac{\partial F_i}{\partial X}\frac{\partial F_j}{\partial Y}\\
\qquad\qquad\qquad\qquad+\frac{\partial F_i}{\partial Y}\frac{\partial F_j}{\partial X})\Gamma[X,Y]dt
+\frac{\partial F_i}{\partial Y}(S_t)\frac{\partial F_j}{\partial Y}(S_t)\Gamma[Y](S_t)dt         \\
\end{array}
\right.
\end{equation}
so that the change of variable formulae for the biases and the variances of errors   (10) are fulfilled if we were defining 
$A[F_i](x,y)$ and $\Gamma[F_i,F_j](x,y)$ by the relations
$$
\begin{array}{c}
A[F_i](S_t)=\frac{(dF_i(S)_t)^\ast}{dt}\\
\\
\Gamma[F_i,F_j](S_t)=\frac{d<F_i(S),F_j(S)>_t}{dt}
\end{array}
$$
the same for  $H$ and any regular function of  $X$ and $Y$.\\

\noindent\underline{To sum up}

We see that the notion of erroneous random variable  $(X,Y)$ with values in  $\mathbb{R}^2$ may represented by 

\noindent a random differential operator of order 1 

$$F\mapsto b_{(X,Y)}[F]=\frac{\partial F}{\partial X}(X,Y)dX+\frac{\partial F}{\partial Y}(X,Y)dY$$

\noindent and a random differential operator of order 2 

$$F\mapsto a_{(X,Y)}[F]=\frac{\partial F}{\partial x}d^2X+\frac{\partial F}{\partial y}d^2Y+\frac{\partial^2 F}{\partial x^2}(dX)^2+2\frac{\partial^2 F}{\partial y\partial y}dXdY
+\frac{\partial^2 F}{\partial y^2}(dY)^2$$
where  $dX,dY,d^2X,d^2Y$ are a priori any random variables.

Then, for an error calculus retaining only the two first moments of the conditional law of the error given   $(X,Y)$, we sum up by two (deterministic) differential operators 

\noindent the bias operator
$$A[F](x,y)=\mathbb{E}[b_{(X,Y)}+\frac{1}{2}a_{(X,Y)}[F]|(X,Y)=(x,y)]$$
and the operator of variance
$$\Gamma[F](x,y)=\mathbb{E}[(b_{(X,Y)}[F])^2|(X,Y)=(x,y)].$$
and we have the change of variable formulae 

\begin{equation}
\label{(10)}
\left\{
\begin{array}{rl}
A[H]=&\sum_{i=1}^d\frac{\partial \Phi}{\partial F_i}A[F_i]+\frac{1}{2}\sum_{i,j=1}^d\frac{\partial^2\Phi}{\partial F_i\partial F_j}\Gamma[F_i,F_j]\\
&\\
\Gamma[H]=&\sum_{i,j=1}^d\frac{\partial \Phi}{\partial F_j}\frac{\partial \Phi}{\partial F_j}\Gamma[F_i,F_j].
\end{array}
\right.
\end{equation}

a) Let us remark that we have identically 
$$(b_{(X,Y)}[F])^2=\frac{1}{2}a_{(X,Y)}[F^2]-Fa_{(X,Y)}[F]$$
therefore
\begin{equation}
\label{ }
\Gamma[F]=A[F^2]-2FA[F]
\end{equation}
the operator $\Gamma$ may be deduced from the operator $a_{(x,y)}$. 

b) If instead of  $a_{(x,y)}$ we had taken $$\hat{a}_{(X,Y)}= c_{(X,Y)}+a_{(X,Y)}$$
where $c$ is a random field of first order differential operators, and so
$$\hat{A}[F]=\mathbb{E}[b_{(X,Y)}[F]+\frac{1}{2}\hat{a}_{(X,Y)}[F]\;|\;(X,Y)\!=\!(x,y)]$$
formulae (12) and (13) would be still verified with  $\hat{A}$, $\Gamma$  remaining unchanged.\\

Remarks a) and b) above show that, for propagating the errors, if we keep the only operators  $A$ and $\Gamma$ and formulae (12), then operator $A$ discloses some   \underline{ambiguity} for its first order terms. We might forget the first order operator $b$. This may be understood by the very nature of the notion of bias which needs an origin for reference. This reference changes if we add at the beginning a deterministic error, of parallax type, which propagates following a first order operator.
\newpage

\noindent{\bf d)} \underline{In manifolds}.

Let us recall that given a manifold  $M$, a \textit{second order tangent vector} at $a$ is a differential operator at point $a$, without constant term, of order $\leq 2$ , we denote $\tau_a(M)$ their set and $\tau(M)$ the space of  fields of second order vectors. 

A second order differential form is a   $\mathcal{C}^\infty$-function on $\tau(M)$ linear on each  $\tau_a(M)$. We denote $\tau^\ast(M)$ the space of second order differential forms. 

Example : If $f$ and $g$ are real functions on $M$ we may define the second order form 
$d^2f=\lambda\in\tau(M)\mapsto\lambda(f)$
and then the second order form $$\qquad df\cdot dg=\frac{1}{2}(d^2(fg)-fd^2g-gd^2f).
$$
We say that $\lambda\in\tau(M)$ is of elliptic type if  $<\!\lambda,df\cdot df\!>\;\geq 0\quad\forall f$. We denote $\tau^e(M)$ the space of fields of elliptic type  second order vectors. \\

We consider a random variable  $X$ with values in  $M$ and, at each point $x\in M$, a probability measure on $\tau_x^e(M)$ regular with respect to $x$.

In other words, we consider a random field of second order tangent vectors of elliptic type (only its marginal laws of order one will be used afterwards). We denote this field $\Delta X$~: it is  the {\tt "}error{\tt "} on $X$.

Then we define
\begin{equation}
\label{variete1}
\left\{\begin{array}{rl}
A[f](x)=&\frac{1}{2}\mathbb{E}[<\!d^2f,\Delta X\!>|X=x]\\
&\\
\Gamma[f](x)=&\mathbb{E}[<\!df\cdot df,\Delta X\!>|X=x]
\end{array}\right.
\end{equation}
 $A$ is a field of second order vectors and $\Gamma$ is a field of bilinear differential operators of positive type given by 
$$\Gamma[f]=A[f^2]-2fA[f].$$
If we consider $h=\varphi(f_1,f_2,\ldots,f_k)$, denoting $x^i$ a coordinate system, we have 
$$\begin{array}{rl}
A[\varphi(f_1,f_2,\ldots,f_k)](x)&=\frac{1}{2}\mathbb{E}_x[<\!\Delta X,d^2h\!>]\\
&=\frac{1}{2}\mathbb{E}_x[<\!\Delta X,\sum_iD_ih\;d^2x^i+\sum_{ij}D_{ij}h\;dx^i\cdot dx^j>]\\
&\\
\Gamma[\varphi(f_1,f_2,\ldots,f_k)](x)&=\mathbb{E}_x[<\!\Delta X,dh\cdot dh\!>]\\
&=\mathbb{E}_x[<\!\Delta X,\sum_{ij}D_ih\;D_jh\;dx"\cdot dx^j\!>]
\end{array}
$$
with
$$\begin{array}{l}
D_ih=\sum_p\varphi^\prime_pD_if_p\\
D_{ij}h=\sum_{p,q}\varphi^{\prime\prime}_{pq}D_if_pD_jf_q+\sum_p\varphi^\prime_pD_{ij}f_p
\end{array}
$$
thus, getting rid :
\begin{equation}
\label{variete2 }
\left\{
\begin{array}{rl}
A[\varphi(f_1,\ldots,f_k)]=&\sum_p\varphi^\prime_pA[f_p]+\frac{1}{2}\sum_{p,q}\varphi^{\prime\prime}_{pq}\Gamma[f_p,f_q]\\
&\\
\Gamma[\varphi((f_1,\ldots,f_k)]=&\sum_{p,q}\varphi^\prime_p\varphi^\prime_q\Gamma[f_p,f_q].
\end{array}\right.
\end{equation}

The interpretation in terms of diffusion processes is similar to that of the flat case. We know that if the Ito differential is   $dF(S)_t$, its finite variation part  $(dF(S)_t)^\ast$ is intrinsic (cf. Paul-Andr\'e Meyer, ``G\'eom\'etrie stochastique sans larmes" {\it in} S\'em de Probabilit\'es XV, Lectures Notes 850, Springer 1981, p 51),  with the interpretation given in the flat case,  $F\mapsto A[F]$ taken in $S_t$ is the second order tangent vector of the local characteristics of  $S$.

Let us remark that if we add to the error  $\Delta X$ a random first order tangent vector, changing $\Delta X$ in  $\Delta X+b$, formula (14) becomes
$$A_b[f](x)=\frac{1}{2}\mathbb{E}_x[<\!df,b\!>+<\!d^2f,\Delta X\!>]$$
because $d^2f|_{TM}=df$ (cf. P.-A. Meyer \textsl{ op. cit.}, p 49) what does not change $\Gamma$ and $A_b$ still satisfies  formulae (15).

Hence we see that the interpretation in terms of errors agrees with the fact that on a manifold, the first order part of a second order tangent vector is not defined [except if a linear connection is available to share the deterministic and the stochastic part (continuous local martingale) of the error]. Already in the flat case, the first order part of the bias needs a convention at the start. \\

\noindent\underline{Remark.} Why did we limit ourselves retaining from the conditional law of the error only the two first moments, since the formalism of $n$-th order differential forms would make it possible to deal with the propagation of the other moments ?
First let us say for the sake of simplicity, but a more precise answer will be given later on ( Remark in part IV page 23).\\

\newpage
{\bf\noindent\textsf{\LARGE III. Complete and symmetric error calculus.}}\\

\noindent a) \underline{The symmetry as invariant.}

If we start from a situation where the error is centered, by a non-linear mapping the error is no more centered. 

What is preserved by image ?

If $X_n$ is an approximation of $X$, and if the joint law of the pair  $(X,X_n)$ is symmetric, then the law of  $(\psi(X),\psi(X_n))$ is of course symmetric too. Concretely this represents situations where we ignore between two close random variables $X_0$ and $X_1$ which is the right one which is the erroneous one, and by this hesitation we work with  $(X_B,X_{1-B})$ where  $B$ is a Bernoulli independent variable.\\

Let us take the notation of the flat case. If $(X,X+\Delta X)$ is a symmetric pair, we have 
$$\begin{array}{c}
F(X+\Delta X)-F(X)=dF+\frac{1}{2}d^2F+\mbox{ remainder }\\
G(X+\Delta X)-G(X)=dG+\frac{1}{2}d^2G+\mbox{ remainder }
\end{array}
$$
so
$$
\begin{array}{l}
\mathbb{E}[F(X\!+\!\Delta X)G(X)\!-\!F(X)G(X\!+\!\Delta X)]\!\\
\qquad\qquad\qquad\qquad=\mathbb{E}[G(X)(dF+\frac{1}{2}d^2F)-F(X)(dF+\frac{1}{2}d^2G)+\mbox{ remainder }]\\
\qquad\qquad\qquad\qquad=\mathbb{E}[G(X)A[F](X)]-\mathbb{E}[F(X)A[G](X)]+\mbox{ remainder }\\
\qquad\qquad\qquad\qquad=<G,AF>_\nu-<F,AG>_\nu+\mbox{ remainder }
\end{array}
$$
denoting $\nu$ the law of $X$.

We see that as soon as the pair $(X,X+\Delta X)$ is symmetric and the errors small, 

\noindent\underline{the operator $A$ is  symmetric with respect to the law of  $X$}. This property of symmetry  (which solves the ambiguity of the first order terms of  $A$) is an invariant by image : the bias operator of the error of the image of $X$ by an application  $\varphi$ is symmetric with respect to the law of  $\varphi(X)$.\\

The symmetric case allows to construct a much more powerful framework thanks to the theory of Dirichlet forms, and this framework extends easily to the infinite dimension  (cf. N. Bouleau, \textsl{Error calculus for finance and physics}, De Gruyter 2003). We sketch it now~:\\
\newpage
\noindent b) \underline{Error structures}

{ An error structure is a term
$$(\Omega, {\cal A}, \mathbb{P}, \mathbb{D}, \Gamma)$$
where $(\Omega, {\cal A}, \mathbb{P})$ is a probability space, satisfying the following properties~:

\noindent1.) $\mathbb{D}$ is a dense sub-vectorspace of   $L^2(\Omega, {\cal A}, \mathbb{P}),$

\noindent2.) $\Gamma$ is a symmetric bilinear map from  $\mathbb{D}\times\mathbb{D\mathbb{}}$ into $L^1(\mathbb{P})$ satisfying
the Gauss calculus of class  ${\cal C}^1\cap {\mbox{Lip}}$, 

\noindent i.e.  if $u\in\mathbb{D}^m$ and $v\in\mathbb{D}^n$,
 and if $F$ and $G$ are of class ${\cal{C}}^1$ and  Lipschitz, 
from
$\mathbb{R}^m$ {\rm[}resp. $\mathbb{R}^n${\rm]} into $\mathbb{R}$, 
 then $F\circ u\in\mathbb{D}$ and $G\circ v\in\mathbb{D}$ and
$$\Gamma[F\circ u,G\circ v]=\sum_{i,j} F_i^{\prime}(u) G_j^{\prime}(v) \Gamma[u_i,v_j]\quad\mathbb{P}{\mbox{-p.s.}},$$
3.) the bilinear form ${\cal E}[f,g]=\mathbb{E}\Gamma[f,g]$ is closed,

\noindent i.e. $\mathbb{D}$ is complete for the norm 
$\|\,.\,\|_{\mathbb{D}}=(\|\,.\,\|_{L^2(\mathbb{P})}^2 +{\cal{E}}[\,.\,,\,.\,])^{\frac{1}{2}}$,

\noindent 4.)  $1\in\mathbb{D}$ (hence $\Gamma[1,1]=0,$ markovianity).

\vspace{.5cm}

\noindent We write ${\cal E}[f]$ for ${\cal E}[f,f]$ and $\Gamma[f]$ for $\Gamma[f,f]$.

\vspace{.5cm}

\noindent With this definition,  the form ${\cal E}$ is a {\it Dirichlet form} (local and  admitting a square field operator). 

To this form corresponds a Dirichlet operator  $A$ (generator of the semi-group associated to   ${\cal E}$) which, under suitable hypotheses on $F$, satisfies
:
$$A[F\circ u]= \sum_i F_i^{\prime}\circ u\;\;A[u_i]+\frac{1}{2}\sum_{i,j}F^{\prime\prime}_{ij}\circ u\;\; \Gamma[u_i,u_j]
\quad\mathbb{P}{\mbox{-p.s.}}.$$}

\newpage
\noindent
{\Large Example 1.}\quad( Ornstein-Uhlenbeck structure in dimension 1)
\begin{eqnarray*}
\Omega & = &\mathbb{R}  \\
\mathcal{A} & = & \hbox{Borel }\sigma{\mbox{-field}} \ \mathcal{B} (
\mathbb{R} )\\
\mathbb{P} & =& \mathcal{N}(0,1) \ \hbox{reduced normal law} \\
\mathbb{D}  = H^1 \bigl( \mathcal{N}(0,1)\bigr) & =& \bigl\{ u\in L^2
({\mathbb{P}} ), u' \ \hbox{in distributions sense}\\
& \ & \qquad\qquad \hbox{belongs to} \ L^2 (\mathbb{P} )\bigr\} \\
\Gamma [u] & = & u^{\prime 2}
\end{eqnarray*}
then, $\bigl( \mathbb{R} ,\mathcal{B} (\mathbb{R} ),
\mathcal{N}(0,1), H^1 (\mathcal{N}(0,1)), \Gamma \bigr)$ is an error structure. We can also obtain the bias operator (the associated generator)~:
$$
\mathcal{D} A = \bigl\{ f\in L^2 (\mathbb{P})\colon \hbox{$f'' -xf'$ in the sense of distributions} \in L^2 (\mathbb{P})\bigr\}
$$
and
$
Af =\frac{1}{2} f'' -\frac{1}{2} I\cdot f'
\quad$
where $I$ is the identity map on $\mathbb{R}$.\\

\noindent
{\Large Example 2.}\quad ( Monte-Carlo structure in dimension 1)
$$\Omega  = [0,1]$$
$$\mathcal{A}  = \hbox{Borel }\sigma{\mbox{-field}}$$
$$\mathbb{P} = \hbox{Lebesgue measure} $$
$$\mathbb{D}  = \bigl\{ u\in L^2 \bigl( [0,1] ,dx\bigr)\colon
\hbox{the derivative $u'$ in the sense of distributions}$$ 
$$\qquad \qquad \hbox{on $]0,1[$ belongs to }
L^2 ([0,1], dx)\bigr\}$$
$$\Gamma [u] = u^{\prime 2}.$$
the space $\mathbb{D}$  thus defined is denoted  $H^1 \bigl( [0,1]\bigr)$.\\

{\Large\noindent Example 3.} Let $U$ be a domain (connected open set) of
$\mathbb{R}^d$ with unit volume,  $\mathcal{B}(U)$ the Borel $\sigma$-field and   $dx =dx_1 ,\ldots dx_d$ the
 Lebesgue measure,
\begin{eqnarray*}
\mathbb{D} & =& \bigl\{ u\in L^2 (U, dx)\colon
\hbox{the gradient $\nabla u$ in the sense of distributions} \\
& \ & \qquad \qquad \hbox{belongs to } L^2 \bigl( U,dx;
\mathbb{R}^d \bigr)\bigr\} \\
\Gamma [u] & =& |\nabla u|^2 = \left( \frac{\partial u}{\partial x_1}
\right)^2 +\cdots +\left( \frac{\partial u}{\partial x_d} \right)^2 .
\end{eqnarray*}
Then $(U, \mathcal{B} (U), dx, \mathbb{D},\Gamma )$ is an error structure.
 From the relation $\mathcal{E} [f,g] =
\langle -Af,g \rangle$ it follows that the domain of the generator contains the functions of class  $\mathcal{C}^2$ with compact support in 
 $U$, $\mathcal{D}A\supset
\mathcal{C}^2_K (U)$ and that for such functions 
$$
Af =\frac{1}{2} \Delta f =\frac{1}{2} \sum^d_{i =1}
\frac{\partial^2 f}{\partial x^2_i} \, .
$$

{\Large\noindent Example 4.}

Let $D$ be an open set in   $\mathbb{R}^d$ with unit volume. Let $\mathbb{P} =dx$ be the Lebesgue measure on $D$.
Let $\Gamma$ be defined on  $\mathcal{C}_K^\infty (D)$ by
$$
\Gamma [u,v] =\sum_{ij} \frac{\partial u}{\partial x_i}
\,\frac{\partial v}{\partial x_j} a_{ij} ,\quad u,v\in \mathcal{C}_K^\infty
(D)
$$
where the functions $a_{ij}$ satisfy the following hypotheses 

\begin{itemize}
\item 
${\displaystyle a_{ij} \in L^2_{\mathrm{loc}} (D) \quad \frac{\partial a_{ij}}{\partial
x_k} \in L^2_{\mathrm{loc}} (D) \quad i,j,k =1,\ldots ,d}$

\item ${\displaystyle \sum_{i,j} a_{ij} (x) \xi_i\xi_j \geq 0
\quad \forall \xi \in \mathbb{R}^d \quad \forall x\in D}$
\item $a_{ij} (x) =a_{ji} (x) \quad \forall x\in D$.
\end{itemize}

\noindent
Then the pre-structure $\bigl( D,\mathcal{B} (D), \mathbb{P} ,
\mathcal{C}_K^\infty (D), \Gamma \bigr)$ is closable.\\

{\sf\Large\noindent Comment on the  completeness. }\\

The fact that we ask the form $\mathcal{E}$ be \textit{closed} is a restriction, but this restriction is highly fruitful.  The situation is quite analogous to the question of the $\sigma$-additivity in probability theory~: without this property, nothing can be said on the transmission of errors by objects defined by limits. But many objects in contemporary mathematics are defined by limits  (integrals, solutions to ode, solutions to sdp, stochastic integrals, solutions to sde, etc.)\footnote{The philosopher Karl Popper is fallen in this trap by emphasizing that his own theory of probability  (additive) do contain strictly the one of Kolmogorov ($\sigma$-additive)
Cf. N.Bouleau
"Some thoughts upon axiomatized languages, a focus on probability theory and error calculus with Dirichlet forms" 
Butlleti de la Societat Catalana de Matem\`atiques Vol. 18 n¡2 p25-36, (2004)
cf. {\tt http://hal.archives-ouvertes.fr/hal-00105636}}

This extension tool allows

- to extend the error calculus on variances from  $\mathcal{C}^1$-functions to  Lipschitz functions,

- to prove a criterion of existence of density for the probability laws which generalizes Malliavin's method~: the image energy density property. 

The closedness property is preserved by images and products even infinite products. This permits to construct naturally errors structures on the spaces of stochastic processes. 

In particular on the Wiener space, where may be obtained, among others, the Ornstein-Uhlenbeck operator. As well, the Malliavin calculus may be interpreted as an error calculus. \\

\noindent
{\bf Theorem on products}

\noindent{\it Let
 $S_n =\bigl( \Omega_n,
\mathcal{A}_n ,\mathbb{P}_n ,\mathbb{D}_n \Gamma_n\bigr)$, $n\geq 1$,
be error structures. The product structure}
$$
S = (\Omega , \mathcal{A}, \mathbb{P} ,\mathbb{D} ,\Gamma )=
\prod^\infty_{n =1} S_n
$$
{\it is defined by}

$$(\Omega ,\mathcal{A} ,\mathbb{P} ) = \left( \prod^\infty_{n =1}
\Omega_{n} , \otimes^\infty_{n =1} \mathcal{A}_n ,
\prod^\infty_{n =1} \mathbb{P}_n\right) $$
\begin{eqnarray*}
\mathbb{D}  =& \biggl\{ f\in L^2 (\mathbb{P} )\colon
\forall n, \ \hbox{\it for almost every } w_1,w_2,\ldots ,w_{n -1},
w_{n +1} ,\ldots\\
 &\hbox{\it for the product measure} \\
 &  x\mapsto f\bigl( w_1 ,\ldots ,w_{n -1} ,x,w_{n +1},
\ldots \bigr) \in \mathbb{D}_n \ \ \hbox{\it and} \\
 &  \int \sum_n \Gamma_n [f] \, d\mathbb{P} < +\infty
\biggr\} 
\end{eqnarray*}
{\it and for $f\in \mathbb{D}$}
$$
\Gamma [f] =\sum^\infty_{n =1}\Gamma_n [f] .
$$
{\it $S$ is an error structure, Markovian if each  $S_n$
is.}\\

\medskip
{\large When we write $\Gamma_n [f]$, $\Gamma_n$ acts only on the 
$n$-th  argument of $f$.}\\

Starting from the basic brics that are one dimensional errors structures,  by product we obtain error structures on function spaces.\\
This yields easily error structures on 

\hspace{1cm} - the Wiener space,   {\rm cf. N. B. \& F. Hirsch, {\sl Dirichlet Forms and Analysis
 on Wiener Space}, De Gruyter 1991},

\hspace{1cm} - the general Poisson space,   the Monte Carlo space,  cf. N. B. \textsl{Error Calculus for Finance and Physics}, De Gruyter, 2003.\\

{\sf\Large\noindent Images of error structures.}\\

The operation is as mere and almost as general as the image of a probability measure by a measurable map. 

If $(\Omega, {\cal A}, \mathbb{P}, \mathbb{D}, \Gamma)$ is an error structure and 
 $X$ a random variable with values in   $\mathbb{R}^d$ whose components are in   
  $\mathbb{D}$, the term $(\mathbb{R}^d, {\cal B}(\mathbb{R}^d),\mathbb{P}_X,\mathbb{D}_X,\Gamma_{\!X})$
is an error structure, where 
$$\begin{array}{rcl}
\mathbb{P}_X&\mbox{is}&\mbox{the law of }X,\\
\mathbb{D}_X&=&\{f\in L^2(\mathbb{P}_X) : f\circ X\in\mathbb{D}\}\\
\Gamma_{\!X}[f](x)&=&\mathbb{E}\{\Gamma[f\circ X]|X\!=\!x\},\quad f\in\mathbb{D}.\end{array}$$
In fact it is possible to define images by more general random variables. \\

The image structure by $X$ may be called the {\tt "}Dirichlet-law{\tt "} of $X$. It is an error structure on $\mathbb{R}^d$ such that the identity map
 $I$ has its components in the domain of  $\Gamma_{\!X}$ and the following formulae hold~:
$$\begin{array}{rl}
\mathbb{E}[\Gamma[X]|X]&= \Gamma_{\!X}[I]\circ X\\
\mathbb{E}[\Gamma[\varphi(X)]|X]&=\Gamma_{\!X}[\varphi]\circ X.
\end{array}$$

Several theorems of probability theory possess analogs in the theory of errors structures under some conditions~:  Gateaux-L\'evy, Strassen, etc. {\rm(cf. Bouleau-Hirsch [de Gruyter 1991].)}

\newpage
{\sf\Large\noindent The case of Wiener space.}\\

\noindent$\bullet$ Let $\bigl( \chi_n\bigr)_{n\in \mathbb{N}}$ be an orthonormal basis of  $L^2 ( \mathbb{R}_+,
\mathcal{B} (\mathbb{R}_+ ), dx)$, and let   $\bigl( g_n\bigr)_{n\in
\mathbb{N}}$ be a sequence of i.i.d. reduced Gaussian random variables.

To a function $f\in L^2 (\mathbb{R}_+,
\mathcal{B} (\mathbb{R}_+ ), dx )$ is associated the Wiener integral 
$$
I(f) =\sum_n \langle f,\chi_n\rangle g_n,
$$
 homomorphism of $L^2 (\mathbb{R}_+, dx)$ into $L^2 (\Omega ,
\mathcal{A} ,\mathbb{P} )$.
If we put
$$
B(t) =\sum_n \bigl\langle 1_{[0,t]} , \chi_n \bigr\rangle
g_n =\sum_n \int^t_0 \chi_n (y)\, dy\cdot g_n
$$
then $B(t)$ is a standard Brownian motion. 

Because of the case where  $f$ is a step-function,
$I(f)$ is denoted
$
I(f) =\int^\infty_0 f(s) \, dB_s 
$\\

\noindent$\bullet$ The preceding construction uses the product space 
$$
(\Omega ,\mathcal{A} ,\mathbb{P} ) =\bigl( \mathbb{R} ,\mathcal{B} 
(\mathbb{R} ), \mathcal{N}(0,1) \bigr)^{\mathbb{N}} ,
$$
the $g_n$'s being the coordinate maps.
If  we put an error structure on each factor
$$
\bigl( \mathbb{R} ,\mathcal{B} (\mathbb{R} ), \mathcal{N}(0,1),
\mathbf{d}_n ,\gamma_n\bigr),
$$
we obtain an error structure on  $(\Omega ,\mathcal{A},
\mathbb{P} )$:
$$
(\Omega ,\mathcal{A} ,\mathbb{P} ,\mathbb{D}, \Gamma )=
\prod^\infty_{n =0} \bigl( \mathbb{R} ,\mathcal{B} (\mathbb{R} ),
\mathcal{N}(0,1) ,\mathbf{d}_n ,\gamma_n \bigr)
$$
such that a random variable 
$
F\bigl( g_0 ,g_1,\ldots ,g_n,\ldots \bigr)
$
belongs to $\mathbb{D}$ iff
$$\left\{
\begin{array}{l}
\forall n\quad x\mapsto F\bigl( g_0,
\ldots ,g_{n -1}, x, g_n,\ldots \bigr) \mbox{ belongs to }\mathbf{d}_n\quad\mathbb{P}\mbox{-p.s.}\\
\mbox{and}
\quad\Gamma [F] =\sum_n \gamma_n [F]\;
\mbox{belongs to }L^1(\mathbb{P})\;(\gamma_n \mbox{ acting on the n-th variable of }F).
\end{array}\right.
$$

\newpage

\noindent$\bullet$ If on each factor we take the Ornstein-Uhlenbeck structure in dimension 1, we obtain 
\begin{eqnarray*}
\Gamma \bigl[ g_n\bigr] & =& 1 \\
\Gamma \bigl[ g_m,g_n\bigr] & =& 0 \quad \hbox{if } m\neq n .
\end{eqnarray*}
For $f\in L^2 (\mathbb{R}_+)$, from $\;\int^\infty_0 f(s)\, dB_s =
\sum\limits_n \langle f,\chi_n\rangle g_n\;$ we deduce
$$
\Gamma \left[ \int^\infty_0 f(s)\, dB_s\right] =\sum_n
\bigl\langle f,\chi_n\bigr\rangle^2 =\| f\|^2_{L^2 (
\mathbb{R}_+)},
$$
hence, using the functional calculus  $\forall F\in \mathcal{C}^1 \cap \hbox{Lip} (\mathbb{R}^m)$
$$\hspace{-.3cm}
\Gamma \left[ F\left( \int f_1 (s)\, dB_s,\ldots ,\int f_n(s)\,
dB_s\right)\right] =\sum_{i,j} \frac{\partial F}{\partial x_i}
\frac{\partial F}{\partial x_j} \int f_i (s)f_j(s)\, ds.
$$
It is the Ornstein-Uhlenbeck structure on the Wiener space. \\

\noindent$\bullet$ For the sake of simplicity let us restrict the time to $t\in[0,1]$ and let us take for $\chi_n$ the trigonometric basis. If on each factor we take the Ornstein-Uhlenbeck structure scaled by a constant coefficient depending on $n$~:
$$
\prod^\infty_{n =0} \bigl( \mathbb{R} ,\mathcal{B} (\mathbb{R} ),
\mathcal{N}(0,1) ,H^1(\mathcal{N}(0,1)) ,u\rightarrow (2\pi n)^{2q}u^{\prime 2} \bigr)
$$
we obtain on the Wiener space an error structure satisfying on the first chaos~:
$$
\Gamma \left[ \int^1_0 f(s)\, dB_s\right]=
\int_0^1 (f^{(q)})^2(s)\,ds
$$
where  $f^{(q)}$ is the $q$-th derivative of $f$. It is a structure where the error disturbs longitudinally the Brownian path, which belongs to the family of structures called of {\tt "}generalized Mehler type{\tt "}. 
\newpage

{\bf\noindent\textsf{\LARGE IV. The four bias operators.}}\\

\noindent We tackle now the following question~: How an error, in the usual sense in mathematics, i.e. an approximation error, generates an error structure. \\

{\large Let us consider a random variable $Y$ defined on the probability space $(\Omega, {\cal A}, \mathbb{P})$  with values in the measurable space 
$(E,{\cal F})$ and approximations $Y_n$, $n\in \mathbb{N}$, also  defined on  $(\Omega, {\cal A}, \mathbb{P})$ with values in $(E,{\cal F})$.

We suppose there exist an algebra ${\cal D}$ of bounded functions from $E$ into $\mathbb{R}$ or $\mathbb{C}$ dense in $L^2(E,{\cal F},\mathbb{P}_Y)$ containing the constants and a sequence  $(\alpha_n)_{n\in\mathbb{N}}$
of positive numbers, about which we consider the following hypotheses~:  
$$
(\mbox{H}1)\qquad\left\{\begin{array}{l}
\forall \varphi\in{\cal D}, \mbox{ there exists } \overline{A}[\varphi]\in L^2(E,{\cal F},\mathbb{P}_Y)\quad \mbox{s.t.} \quad\forall \chi\in{\cal D}\\
\lim_{n\rightarrow\infty} \alpha_n\mathbb{E}[(\varphi(Y_n)-\varphi(Y))\chi(Y)]=\mathbb{E}_Y[\overline{A}[\varphi]\chi]
\end{array}\right.
$$
the expectation $\mathbb{E}_Y$ being relative to the law $\mathbb{P}_Y$.
$$
(\mbox{H}2)\qquad\left\{\begin{array}{l}
\forall \varphi\in{\cal D}, \mbox{ there exists } \underline{A}[\varphi]\in L^2(E,{\cal F},\mathbb{P}_Y)\quad \mbox{s.t.} \quad\forall \chi\in{\cal D}\\
\lim_{n\rightarrow\infty} \alpha_n\mathbb{E}[(\varphi(Y)-\varphi(Y_n))\chi(Y_n)]=\mathbb{E}_Y[\underline{A}[\varphi]\chi].
\end{array}\right.
$$
$$
(\mbox{H}3)\quad\left\{\begin{array}{l}
\forall \varphi\in{\cal D}, \mbox{ there exists } \widetilde{A}[\varphi]\in L^2(E,{\cal F},\mathbb{P}_Y)\quad \mbox{s.t.} \quad\forall \chi\in{\cal D}\\
\lim_{n\rightarrow\infty} \alpha_n\mathbb{E}[(\varphi(Y_n)-\varphi(Y))(\chi(Y_n)-\chi(Y))]=-2\mathbb{E}_Y[\widetilde{A}[\varphi]\chi].
\end{array}\right.
$$

As soon as two hypotheses  among (H1) (H2) (H3) are fulfilled  (with the same algebra ${\cal D}$ and the same sequence $\alpha_n$), the third one follows thanks to the relation $$\widetilde{A}=\frac{\overline{A}+\underline{A}}{2}.$$

$\bullet$ When defined the operator $\overline{A}$ which considers the asymptotic error from the point of view of the limit model, will be called  
{\it the theoretical bias operator}.

$\bullet$ The operator  $\underline{A}$ which consider the error from the point of view of the approximate model, will be called  
 {\it the practical bias operator}.

$\bullet$ Because of the property
$$<\widetilde{A}[\varphi],\chi>_{L^2(\mathbb{P}_Y)}=<\varphi,\widetilde{A}[\chi]>_{L^2(\mathbb{P}_Y)}$$
the operator $\widetilde{A}$ will be called {\it the symmetric bias operator}.\\

The following result shows that a Dirichlet form (possibly non-local) is often present behind an approximation~:\\

{\large\noindent{\bf Theorem. }{\it  Under hypothesis} (H3) {\it

a) the limit
$$\widetilde{\cal E}[\varphi,\chi]=\lim_n  \frac{\alpha_n}{2}\mathbb{E}[(\varphi(Y_n)-\varphi(Y))(\chi(Y_n)-\chi(Y)]\qquad \varphi, \chi\in{\cal D}$$
defines a closable bilinear form whose smallest closed extension is denoted  $({\cal E},\mathbb{D})$.

b) $({\cal E},\mathbb{D})$ is a Dirichlet form

c) $({\cal E},\mathbb{D})$ admits a square of field operator  $\Gamma$ satisfying $\forall \varphi,\chi\in{\cal D}$
$$
\Gamma[\varphi]=\widetilde{A}[\varphi^2]-2\varphi\widetilde{A}[\varphi]$$
$$\mathbb{E}_Y[\Gamma[\varphi]\chi]=\lim_n\alpha_n\mathbb{E}[(\varphi(Y_n)-\varphi(Y))^2(\chi(Y_n)+\chi(Y))/2]$$
\indent d) $({\cal E},\mathbb{D})$ is local if and only if  $\forall \varphi\in{\cal D}$
$$\lim_n \alpha_n\mathbb{E}[(\varphi(Y_n)-\varphi(Y))^4]=0.$$}}

We introduce now the fourth bias operator  $\Abar$ defined under  (H1) and (H2) on ${\cal D}$ by
$$\Abar=\frac{1}{2}(\overline{A}-\underline{A}).$$
Since $\mathbb{E}_Y[\Abar[\varphi]\chi]=\lim_n\mathbb{E}[(\varphi(Y_n)-\varphi(Y))(\chi(Y)+\chi(Y_n))/2]$ we see that $\Abar$ represents the asymptotic error from the point of view of an outside observer  who  attaches the same weight both to the theoretical and the practical models and measuring the error algebraically on the same axis.  Because of properties we shall state below,  the operator 
 $\Abar$ will be called
  {\large\it the singular bias operator}.\\

An operator  $B$ from ${\cal D}$ into $L^2(\mathbb{P}_Y)$ is said {\it first order} if it satisfies
$$B[\varphi\chi]=B[\varphi]\chi+\varphi B[\chi]\qquad\forall\varphi,\chi\in{\cal D}.$$
\newpage
\noindent{\bf Proposition.} {\it Under} (H1) {\it to} (H3) {\it 

a) the}  theoretical variance $\lim_n\alpha_n\mathbb{E}[(\varphi(Y_n)-\varphi(Y))^2\psi(Y)]$ {\it and the } practical variance
$\lim_n\alpha_n\mathbb{E}[(\varphi(Y_n)-\varphi(Y))^2\psi(Y_n)]$ {\it exist and we have $\forall\varphi,\chi,\psi\in{\cal D}$
$$
\begin{array}{c}
\lim_n\alpha_n\mathbb{E}[(\varphi(Y_n)-\varphi(Y))(\chi(Y_n)-\chi(Y))\psi(Y)]=\mathbb{E}_Y[-\underline{A}[\varphi\psi]\chi+\underline{A}[\psi]\varphi\chi
-\overline{A}[\varphi]\chi\psi]\\
\lim_n\alpha_n\mathbb{E}[(\varphi(Y_n)-\varphi(Y))(\chi(Y_n)-\chi(Y))\psi(Y_n)]=\mathbb{E}_Y[-\overline{A}[\varphi\psi]\chi+\overline{A}[\psi]\varphi\chi
-\underline{A}[\varphi]\chi\psi]
\end{array}
$$
\indent \it b) These two variances  coincide if and only if $\;\Abar$ is first order, and then they are equal to 
$\mathbb{E}_Y[\Gamma[\varphi]\psi].$

c) If the Dirichlet form is local, then $\;\Abar$ is first order.
}\\

\noindent\underline{Remark.} Under (H3) condition d) of the theorem (p 22)  caracterising the case where the form ${\cal E}$ is local is equivalent to either one of the following conditions~: 

(j) $\exists\lambda>2\quad\lim_n\alpha_n\mathbb{E}[|\varphi(Y_n)-\varphi(Y)|^\lambda]=0\quad\forall\varphi\in{\cal D}.$

(jj) $\forall\lambda>2\quad\lim_n\alpha_n\mathbb{E}[|\varphi(Y_n)-\varphi(Y)|^\lambda]=0\quad\forall\varphi\in{\cal D}.$

\noindent This yields an answer to the question put on page 12 about an error calculus with higher moments using tangent vectors of order greater than 2 ~: In the cases where the Dirichlet form is local, (the only case in which we are able to propagate errors thanks to a differential calculus) the moments of order greater than 2 are negligible before the first or second order moments. \\

Let us come back to the case where only hypothesis (H3) is assumed. The following result shows that, for the variances, the error calculus performed on the limit model with the asymptotic error, coincides with the error asymptotically obtained on  $\mathcal{C}^1$-functions. \\

\noindent{\bf Proposition.} {\it Under} (H3). {\it If the form $({\cal E},\mathbb{D})$ is local, then the} asymptotic error calculus principle 
{\it is valid on
$$\widetilde{\cal D}=\{F(f_1,\ldots,f_p)\;:\;f_i\in{\cal D},\;\;F\in{\cal C}^1(\mathbb{R}^p,\mathbb{R})\}$$
i.e.}
$
\lim_n\alpha_n\mathbb{E}[(F(f_1(Y_n),\ldots,f_p(Y_n))-F(f_1(Y),\ldots,f_p(Y))^2]$\hfill

\hfill$=\mathbb{E}_Y[\sum_{i,j=1}^p F^\prime_i(f_1,\ldots,f_p)F^\prime_j(f_1,\ldots,f_p)\Gamma[f_i,f_j]].$\\

\noindent{\sf Examples}\

$\bullet$ Let us take for $Y$ a Brownian motion $B$ indexed by $[0,1]$ as random variable with value in  ${\cal C}([0,1])$ and let us take for  $Y_\varepsilon$  the approximation $Y_\varepsilon=B+\sqrt{\varepsilon}W$ where $W$ is an independent Brownian motion. We may apply the theorem taking for  ${\cal D}$ the linear combinations  of functions   $\varphi(B)=e^{i\int_0^1f\,dB}$ with regular $f$ say ${\cal C}^1_b$.

Hypothesis (H3) is fulfilled. The theorem yields the Ornstein-Uhlenbeck structure on the Wiener space. 
 
 $\bullet$ Series with independent increments. let be $$S=\sum_{n=1}^\infty \frac{X_n}{n^2}+\frac{Z_n}{n}$$ where $X_n,Z_n\in L^{2+\varepsilon}$, $Z_n$ centered, $(X_n,Z_n)$ i.i.d., on 
and $S$ is approximated by the partial sum  $S_n=\sum_{k=1}^n\frac{X_k}{k^2}+\frac{Z_k}{k}$. \\

By Burkholder's inequality, we have $n\mathbb{E}[|S-S_n|^{2+\varepsilon}]\rightarrow 0$ as $n\rightarrow\infty$. So that, taking ${\cal D}={\cal C}^\infty_K$, for $\varphi,\chi\in{\cal D}$ 
$$\lim_n n\mathbb{E}[(\varphi(S)-\varphi(S_n))\chi(S_n)]=\lim_n n\mathbb{E}[(S-S_n)\varphi^\prime(S_n)\chi(S_n)+\frac{1}{2}(S-S_n)^2\varphi^{\prime\prime}(S_n)\chi(S_n)]$$
$$=\frac{1}{2}\mathbb{E}[Z_1^2]\mathbb{E}[\varphi^{\prime\prime}(S)\chi(S)]+\mathbb{E}[X_1]\mathbb{E}[\varphi^\prime(S)\chi(S)].$$
Hence hypothesis  (H2) is satisfied and $$\underline{A}[\varphi]=
\frac{\mathbb{E}[Z_1^2]}{2}\varphi^{\prime\prime}+ \mathbb{E}[X_1]\varphi^\prime.$$
Similarly
$$\lim_n n\mathbb{E}[(\varphi(S)-\varphi(S_n))^2]=\lim_n n\mathbb{E}[(S-S_n)^2\varphi^{\prime 2}(S_n)]=\mathbb{E}[Z_1^2]\mathbb{E}[\varphi^{\prime 2}(S)]$$ 
 (H3) is fulfilled as soon as the law of  $S$ satisfies the Hamza condition (cf. M. Fukushima et al. [1994]) and then the Dirichlet form is local, and $\Gamma[\varphi]=\underline{A}[\varphi^2]-2\varphi\underline{A}[\varphi]=\mathbb{E}[Z_1^2]\varphi^{\prime2}$.\\
 
 $\bullet$  Stochastic intergral. Let us consider the integral
$$Y=\int_0^1 H_s\,dB_s$$ approximated by the sum 
$$Y_n=\sum_{k=0}^{n-1}H_{\frac{k}{n}}(B_{\frac{k+1}{n}}-B_{\frac{k}{n}})$$
$(B_t)$ is a standard Brownian motion defined as coordinate process on the space  ${\cal C}([0,1])$ equipped with Wiener measure, $H_s=H_0+\int_0^s\xi_u\,dB_u+\int_0^s\eta_u\,du$ is an Ito process regular in Malliavin's sense. Under suitable hypotheses we obtain~:

- Hypothesis (H3) is fullfilled under rather simple regularity assumptions and   $$<\widetilde{A}[\varphi],\chi>=-\frac{1}{4}\mathbb{E}[\int_0^1 \xi_s^2\,ds\varphi^\prime(Y)\chi^\prime(Y)].$$

- Hypothesis (H1) supposes finer regularity conditions and  
$$\begin{array}{rl}
<\overline{A}[\varphi],\chi>&=\frac{1}{2}\mathbb{E}[\int_0^1\xi_sD_sD_s[\varphi^\prime(Y)\chi(Y)]ds]+\frac{1}{2}\mathbb{E}[\int_0^1\eta_sD_s[\varphi^\prime(Y)\chi(Y)]ds]\\
&\quad \frac{1}{4}\mathbb{E}[\int_0^1\xi_s^2(\varphi^\prime\chi)^\prime(Y) ds]-\frac{1}{4}\mathbb{E}[\int_0^1\xi_s^2\varphi^\prime(Y)\chi^\prime(Y) ds]
\end{array}
$$ where $D$ denotes Malliavin's gradient on the Wiener space. Then  $\Abar$ is a first order operator.\\

$\bullet$ Stochastic differential equations and the Euler scheme. The approximation error of the solution of an SDE by the Euler scheme has been the subject of many works whose one of the achievements is a form of functional central limit theorem. See especially { Jean Jacod and Philip Protter} ``Asymptotic error distributions for the Euler method for stochastic 
differential equations'' {\sl  Ann. Probab.} 26, 267-307, (1998) and the quoted references. This central limit theorem allows to simplify the study of the limit
$$\lim_{n\rightarrow\infty} \alpha_n\mathbb{E}[(\varphi(Y_n)-\varphi(Y))^2]$$ but the existence of the operator $\widetilde{A}$, and therefore the closability of the form are yet only proved in dimension 1. On the other hand the operator $\overline{A}$ has been determined by Paul Malliavin and Anton Thalmaier ``Numerical error for SDE: asymptotic expansion and hyperdistributions" {\it Note C. R. A. S.} sI, vol 336, n¡10, p851, 2003, its appearance is similar to the one given above in the case of a stochastic integral. \\

\newpage
{\bf\noindent\textsf{\LARGE V. Statistics and errors.}}\\

In the error calculus by Dirichlet forms  {\it any erroneous quantity is random}. The a priori probability law may be thought as  {\tt "}the field and the use{\tt "} of the measurement device.\\

{\noindent\large How to determine by experiment the operator  $\Gamma$ of an error structure ?}\\

 Let us suppose that the error structure to be identified be on $\mathbb{R}^d$:
$$(\mathbb{R}^d, {\cal B}(\mathbb{R}^d), \mathbb{P}, \mathbb{D}, \Gamma).$$
Concretely  $(\mathbb{R}^d, {\cal B}(\mathbb{R}^d), \mathbb{P})$ is the image space of a quantity $x$ which is measured with some accuracy.  $\mathbb{P}$ is its a priori law.\\

We shall consider that to carry out a measurement of the quantity  $x$, is to estimate $x$ statistically as parameter of a family of probability measures 
 $\mathbb{Q}_x$.

It is known that if $T$ is an estimate of  $x$ unbiased say 
\noindent($\mathbb{Q}_x[T]=x$), the precision on $x$
is limited by the   inequality
$$\mathbb{Q}_x[(T-x)(T-x)^t]\geq[J(x)]^{-1}$$ with equality if $T$ is efficient (Cramer-Darmois-Fisher-Rao inequality).
Let us recall the definition of the Fisher information and the Cramer-{\it et-al.} inequality. 


Let be $x\in\mathbb{R}^d$ and let $\mathbb{Q}_x$ be a family of probability measures on some space dominated by the probability measure 
$\mathbb{Q}$
$$\mathbb{Q}_x=L(x,.)\,\mathbb{Q}\quad\mbox{ with } L(x,.) \mbox{ regular in }x.$$
Then for every random variable  $Y\in L^2(\mathbb{Q})$ 
{\large$$\mathbb{E}_x[Y-\mathbb{E}_x(Y)]^2\geq (\mathrm{grad}\mathbb{E}_x(Y))^t[J(x)]^{-1}\mathrm{grad}\mathbb{E}_x(Y)$$}
where $J(x)$ is the Fisher information matrix of the model {\large$$J(x)=\left(\mathbb{E}_x[\frac{\partial\log L(x)}{\partial x_i}\frac{\partial\log L(x)}{\partial x_j}]\right)_{ij}.$$}
Let $I$ be the identity map from $\mathbb{R}^d$ into itself, $J(x)$ behaves as an information,
$\Gamma[I](x)$ is a precision.

We assume the identification~:

{\Large\begin{center}
\begin{tabular}{|c|}
\hline
$\Gamma[I](x)=J^{-1}(x)$\\
\hline
\end{tabular}
\end{center}}


If we choose  this way $\Gamma[I](x)=J^{-1}(x)$, since $\Gamma$ satisfies the functional calculus, if
$f:\mathbb{R}^d\longrightarrow\mathbb{R}^p$ is of class ${\cal C}^1\cap\mathrm{Lip}$, we have
$$\Gamma[f](x) =(\mathrm{grad}f)^t.\Gamma[I](x).\mathrm{grad}f$$
In other words the operator $\Gamma$ is uniquely determined, and similarly for all image structures of our error structure 
$(\mathbb{R}^d, {\cal B}(\mathbb{R}^d), \mathbb{P}, \mathbb{D}, \Gamma)$.\\

Hence the question arises to know whether we obtain a precision compatible with this calculus when   for $f$ one to one we measure $f(x)$.
The answer is poisitive. 

{ Under the hypotheses of a {\tt "}regular{\tt "} statistical model, if we consider $y=f(x)$ for $f\in{\cal C}^1\cap\mathrm{Lip}$
and one to one, the error structure obtained for $y$ is the image of the error structure obtained for  $x$ by the
 map $f$. The identification is also stable by products in a natural sense.} Cf. { N. B. and Chr. Chorro, ``Structures d'erreur et estimation param\`etrique'',
 {\it Note C.R.A.S.}, Ser.I 338 (2004), 305-310}.

That means that the error on $x$ obtained this way {\it do not depend on parametrization} and
{\it possesses a physical meaning}.\\

\noindent {\sf Some books related to these subjects}\\

{\normalsize
{\sc Azema, J.; Yor, M.} {\it S\'eminaire de Probabilit\'es XVI, 1980/81 Suppl\'ement : G\'eom\'etrie Diff\'erentielle Stochastique}, LNM 921, Springer, 1982.

{\sc Bouleau N.} {\it Error Calculus for Finance and Physics, the Language of Dirichlet Forms}, De Gruyter, 2003.

{\sc Bouleau N., Hirsch F.}  {\it Dirichlet Forms and Analysis on Wiener Space,} De Gruyter, (1991).

 {\sc Dellacherie, Cl., Meyer, P.-A.,} {\it Probabilit\'es et Potentiel,} Hermann 1987.
 
 {\sc Emery, M.} {\it Stochastic Calculus in Manifolds},  Springer, 1989.

{\sc Fukushima, M.; Oshima, Y.; Takeda, M.} {\it Dirichlet Forms and Symmetric Markov Processes}, De Gruyter 1994.

{\sc Jacod, J., Shiryaev, A.N.,} {\it Limit Theorems for Stochastic Processes,} Springer 1987.

{\sc Ma, Z.-M., R\"{o}ckner, M.} {\it Introduction to the Theory of (Non-symmetric) Dirichlet Forms,} Springer 1992.

{\sc Malliavin, P., Thalmaier, A.} {\it Stochastic Calculus of Variations in Mathematical Finance}, Springer 2005.

{\sc Nualart, N.} : {\it The Malliavin Calculus and Related Topics}, Springer, 1995.}

\end{document}